\documentclass{conm-p-l}

\usepackage{graphicx}
\usepackage[T1]{fontenc}
\usepackage{tikz}

\begin{document}

\newcommand{\EL}{{ L}}
\newcommand{\R}{\mathbb{R}}
\newcommand{\Q}{\mathbb{Q}}
\newcommand{\Z}{\mathbb{Z}}
\newcommand{\homol}{\widetilde{H}}
\newcommand{\prob}{\mathbb{P}}
\newcommand{\expect}{\mathbb{E}}
\newcommand{\var}{\mbox{Var}}
\newcommand{\cov}{\mbox{Cov}}
\newcommand{\lk}{\mbox{lk}}
\newcommand{\Pois}{\mbox{Pois}}
\newcommand{\link}{\mbox{lk}}
\newcommand{\cd}{\operatorname{cd}}

\newtheorem{theorem}{Theorem}[section]
\newtheorem{lemma}[theorem]{Lemma}
\newtheorem{corollary}[theorem]{Corollary}
\newtheorem{conjecture}[theorem]{Conjecture}
\newtheorem{proposition}[theorem]{Proposition}
\theoremstyle{definition}
\newtheorem{definition}[theorem]{Definition}
\newtheorem{example}[theorem]{Example}
\newtheorem{claim}[theorem]{Claim}

\title{Topology of random simplicial complexes: a survey}


\author{Matthew Kahle}
\address{The Ohio State University, Department of Mathematics}
\email{mkahle@math.osu.edu}
\thanks{The author gratefully acknowledges support from DARPA grant \# N66001-12-1-4226}


\date{\today}

%
%

\maketitle

\section{Introduction}

This expository article is based on a lecture from the Stanford Symposium on Algebraic Topology: Application and New Directions, held in honor of Gunnar Carlsson, Ralph Cohen, and Ib Madsen in July 2012.

\bigskip

\begin{quote}
``I predict a new subject of statistical topology. Rather than count the number of holes, Betti numbers, etc., one
will be more interested in the distribution of such objects on noncompact manifolds as one goes out to infinity.''\\
  --- Isadore Singer, 2004 \cite{interview}
\end{quote}

\bigskip

\subsection{Motivation}

\subsubsection{Randomness models the natural world}

\begin{enumerate}
\item The field of topological data analysis has received a lot of attention over the past several years --- see, for example, the survey articles of Carlsson \cite{C09} and Ghrist \cite{Ghrist}.  In order to quantify the statistical significance of topological features detected with these methods, it will be helpful to have firmer probabilistic foundations.

\item Certain areas of physics seem to require random topology. John Wheeler suggested in the 1960's, for example that inside a black hole, one would need a topological and geometric theory to account for relativity, and that near the presumed singularity one would also require a quantum theory, necessarily stochastic.  He reasoned that inside a black hole the topology of space-time itself may best be described as a quantum foam, or as a probability distribution over shapes rather than any particular fixed shape.

\item We might also want to understand better why certain mathematical phenomena are so ubiquitous.  For example, a folklore theorem is that ``almost all groups are hyperbolic''.  This turns out to be true under a variety of different measures --- see, for example, Ollivier's survey on random groups \cite{Ollivier}.  

It is known that many simplicial complexes and posets found in the wild are homotopy equivalent to bouquets of spheres. See for example Forman's comments in Section 1 of his discrete Morse theory notes \cite{Forman}. He points out that Morse theory gives a convenient way to prove such theorems, but goes on to say, ``However, that does not explain why so many simplicial complexes that arise in combinatorics are homotopy equivalent to a wedge of spheres. I have often wondered if perhaps there is some deeper explanation for this.'' From the point of view of this article, one might hope to make mathematical sense of such questions measure-theoretically.

\item Random topology might even provide tractable toy models for difficult to understand number-theoretic settings.

The following family of simplicial complexes was apparently introduced and first studied topologically by Bj\"{o}rner \cite{B11}:  $\Delta_n$ has primes less than $n$ as its vertices, and its faces correspond to square-free numbers $i$ with $1 \le i \le n$.  He pointed out that the Euler characteristic $\chi ( \Delta_n)$ coincides with the Mertens function $$M(n) = \sum_{k=1}^n \mu(k),$$ where $\mu(k)$ is the M\"{o}bius function.  The Riemann hypothesis is equivalent to the statement that $M(n)$ satisfies
$$| M(n) | = O(n ^{1/2 + \epsilon}) $$
for every $\epsilon > 0$, suggesting that studying the topology of $\Delta_n$ might be quite interesting.

Bj\"{o}rner proved these complexes are all homotopy equivalent to wedges of spheres (not necessarily of the same dimension), hence homology $H_*(\Delta_n)$ is torsion-free.  He also provided some estimates for Betti numbers, showing that $$\beta_k \approx \frac{n}{2 \log n} \frac{(\log \log n)^k}{k!} $$ for $k \ge 0$ fixed, and also that
$$\sum_{k \ge 0} \beta_k \left( \Delta_n \right) = \frac{2n}{\pi^2} + O \left( n^{\theta} \right)$$ for all $\theta > 17/54$.

Pakianathan and Winfree \cite{PW13} studied a fairly general framework ``quota complexes'', and gave natural topological formulations of the prime number theorem, twin prime conjecture, Goldbach's conjecture, and the existence of odd perfect numbers, among others.

Unfortunately, these attractive papers do not seem to get us any closer to proving the Riemann hypothesis, but these kinds of complexes merit further study, and might be interesting to model probabilistically. Of course primes are not random, but they are pseudorandom, and for many purposes behave as if they were a random subset of the integers with density predicted by the prime number theorem --- the Green--Tao theorem is a celebrated example \cite{GT08}.

\end{enumerate}

\subsubsection{The probabilistic method provides existence proofs.}

This is a complementary point of view.  Random objects often have desirable properties, and in some cases it is difficult to construct explicit examples.  This has been one of the most influential ideas in discrete mathematics of the past several decades --- for a broad overview of the subject, see Alon and Spencer's book \cite{Alon}.

\begin{enumerate}
\item In extremal graph theory, for example, the probabilistic method has proved to be an extremely powerful tool.  Almost all graphs are known to have strong Ramsey properties (i.e.\ no large cliques or independent sets), but after several decades of research no one knows how to construct examples.  This somewhat paradoxical situation is sometimes referred to as the problem of ``finding hay in a haystack.'' 

\item Since early work of Pinsker \cite{P73}, and even earlier work of Barzdin and Kolmogorov \cite{BK67}, it has been known that in various senses almost all graphs are expanders.  For an extensive survey of expander graphs and their many applications in mathematics and computer science, see \cite{HLW06}.  Some of the work surveyed in this article may be viewed as higher-dimensional analogues of this paradigm. One of the goals of this article is to describe expander-like qualities of random simplicial complexes.

\item The probabilistic method has found applications in other areas of mathematics as well.  Gromov asked, ``What does a random group look like? As we shall see the answer is most satisfactory: nothing like we have ever seen before,'' \cite{spaces}, and then later fulfilled his own prediction by proving the existence of a finitely generated group $\Gamma$ whose Cayley graph admits no uniform embedding in the Hilbert space \cite{rwirg}.  Rather than exhibit such a group explicitly, he makes a very delicate measure-theoretic argument.
\end{enumerate}

\subsection{Earlier work}

Although this article will focus on random simplicial complexes, we first put this in a larger context of random topology.

\begin{enumerate}

\item One of the earliest references to random topology of which I am aware is Milnor's 1964 paper ``Most knots are wild'' \cite{wild}. Milnor showed that ``most'' knots are wild, in the sense of wild embeddings being Baire dense in the space of all embeddings. In the concluding remarks he points out that questions about embeddings being wild or knotted with probability one are of a very different kind. In particular, he noted that with probability one, Brownian motion in $4$-space does not self intersect and asked if is knotted. Kendall \cite{Kendall85} and independently Weinberger \cite{Weinberger13} showed that the answer is ``no.''

\item Random triangulated surfaces were studied by Pippenger and Schleich \cite{PS}.  Their model is randomly gluing together $n$ oriented triangles, uniformly over all such glueings, and they compute the expected genus $\expect[g_n]$ of the resulting oriented surface as $n \to \infty$.  Part of the motivation discussed comes from physics --- such random surfaces arise in 2-dimensional quantum gravity and as world-sheets in string theory.

\item Dunfield and W.\ Thurston \cite{DT} simultaneously studied this model for random surfaces, and they pointed out that in general one can not make a random $3$-manifold by gluing together tetrahedra in an analogous way, as the probability that a gluing results in a manifold tends to $0$ as the number of tetrahedra $n \to \infty$. They introduced a new model for random $3$-manifolds $M$ where one takes a random walk on the mapping class group, resulting in a random gluing of a two handlebodies. 

\item Farber and Kappeler \cite{FK08} introduced random planar linkages, were the lengths of the links are random. The configuration space for such a linkage is a smooth manifold with probability one. They let the number of links (and hence the dimension of the manifold) tend to infinity, and give asymptotic formula for the expectations of the Betti numbers.

\item Gaussian random functions on manifolds were studied by Adler and Taylor, and in particular they discovered the kinematic formula for the expected Euler characteristic of the sub-level sets.  See for example, Chapter 12 of their book \cite{AT07}.  In this setting, giving a formula for the expectation of the Betti numbers seems to be a difficult open problem.  For a survey of this area see \cite{ABBSW}.

\item Random $2$-dimensional simplicial complexes were first studied by Linial and Meshulam \cite{LM}, and the $k$-dimensional version by Meshulam and Wallach \cite{MW}.  These are natural higher-dimensional analogues of random graphs --- the main results of \cite{LM}, \cite{MW}, and \cite{flag} are cohomological analogues of the Erd\H{o}s--R\'enyi theorem which characterizes the threshold for connectivity of a random graph. All of these theorems describe sharp topological phase transitions where cohomology passes to vanishing with high probability, within a very short window of parameter.

\end{enumerate}

Since the influential papers \cite{LM} and \cite{MW}, random complexes and their topological properties have continued to be explored by several teams of researchers. Random flag complexes  \cite{clique, flag} are another generalization of Erd\H{o}s--R\'enyi random graphs to higher dimensions, putting a measure on a wide range of possible topologies.

\section{Random graphs} \label{sec:rg}

Random graphs are $1$-dimensional random simplicial complexes. The random graph $G(n,p)$, sometimes called the Erd\H{o}s--R\'enyi model, has vertex set $[n] = \{ 1, 2, \dots, n\}$ and every possible edge appears independently with probability $p$. The closely related random graph $G(n,m)$ is selected uniformly among all ${n \choose 2} \choose m$ graphs on $n$ vertices with exactly $m$ edges.

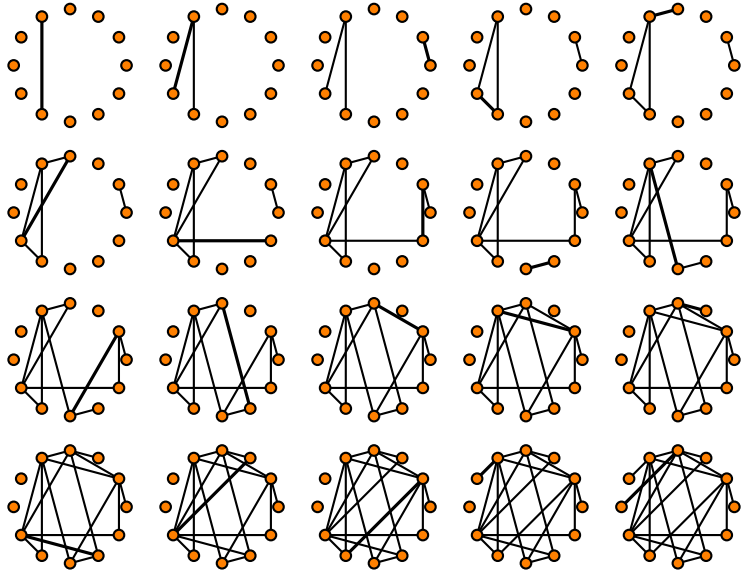
\begin{figure}
\tikzstyle{every node}=[circle, draw, fill=orange,
                        inner sep=0pt, minimum width=4pt]
\begin{tikzpicture}[thick,scale=0.25]
    \draw[very thick] (120:3) -- (240:3);
    \draw \foreach \x in {0,30,...,330}
    {
      (\x:3)  node {}
    };
\end{tikzpicture}\quad
\begin{tikzpicture}[thick,scale=0.25]
    \draw (120:3) -- (240:3);
    \draw[very thick] (120:3) -- (210:3);
    \draw \foreach \x in {0,30,...,330}
    {
      (\x:3)  node {}
    };
\end{tikzpicture}\quad
\begin{tikzpicture}[thick,scale=0.25]
    \draw (120:3) -- (240:3);
    \draw (120:3) -- (210:3);
    \draw[very thick] (0:3) -- (30:3);
    \draw \foreach \x in {0,30,...,330}
    {
      (\x:3)  node {}
    };
\end{tikzpicture}\quad
\begin{tikzpicture}[thick,scale=0.25]
    \draw (120:3) -- (240:3);
    \draw (120:3) -- (210:3);
    \draw (0:3) -- (30:3);
    \draw[very thick] (210:3) -- (240:3);
    \draw \foreach \x in {0,30,...,330}
    {
      (\x:3)  node {}
    };
\end{tikzpicture}\quad
\begin{tikzpicture}[thick,scale=0.25]
    \draw (120:3) -- (240:3);
    \draw (120:3) -- (210:3);
    \draw (0:3) -- (30:3);
    \draw (210:3) -- (240:3);
    \draw[very thick] (90:3) -- (120:3);
    \draw \foreach \x in {0,30,...,330}
    {
      (\x:3)  node {}
    };
\end{tikzpicture}\quad
\\
\vspace{0.1in}
\begin{tikzpicture}[thick,scale=0.25]
    \draw (120:3) -- (240:3);
    \draw (120:3) -- (210:3);
    \draw (0:3) -- (30:3);
    \draw (210:3) -- (240:3);
    \draw (90:3) -- (120:3);
    \draw[very thick] (90:3) -- (210:3);
    \draw \foreach \x in {0,30,...,330}
    {
      (\x:3)  node {}
    };
\end{tikzpicture}\quad
\begin{tikzpicture}[thick,scale=0.25]
    \draw (120:3) -- (240:3);
    \draw (120:3) -- (210:3);
    \draw (0:3) -- (30:3);
    \draw (210:3) -- (240:3);
    \draw (90:3) -- (120:3);
    \draw (90:3) -- (210:3);
    \draw[very thick] (210:3) -- (330:3);
    \draw \foreach \x in {0,30,...,330}
    {
      (\x:3)  node {}
    };
\end{tikzpicture}\quad
\begin{tikzpicture}[thick,scale=0.25]
    \draw (120:3) -- (240:3);
    \draw (120:3) -- (210:3);
    \draw (0:3) -- (30:3);
    \draw (210:3) -- (240:3);
    \draw (90:3) -- (120:3);
    \draw (90:3) -- (210:3);
    \draw (210:3) -- (330:3);
    \draw[very thick] (330:3) -- (30:3);
    \draw \foreach \x in {0,30,...,330}
    {
      (\x:3)  node {}
    };
\end{tikzpicture}\quad
\begin{tikzpicture}[thick,scale=0.25]
    \draw (120:3) -- (240:3);
    \draw (120:3) -- (210:3);
    \draw (0:3) -- (30:3);
    \draw (210:3) -- (240:3);
    \draw (90:3) -- (120:3);
    \draw (90:3) -- (210:3);
    \draw (210:3) -- (330:3);
    \draw (330:3) -- (30:3);
    \draw[very thick] (270:3) -- (300:3);
    \draw \foreach \x in {0,30,...,330}
    {
      (\x:3)  node {}
    };
\end{tikzpicture}\quad
\begin{tikzpicture}[thick,scale=0.25]
    \draw (120:3) -- (240:3);
    \draw (120:3) -- (210:3);
    \draw (0:3) -- (30:3);
    \draw (210:3) -- (240:3);
    \draw (90:3) -- (120:3);
    \draw (90:3) -- (210:3);
    \draw (210:3) -- (330:3);
    \draw (330:3) -- (30:3);
    \draw (270:3) -- (300:3);
    \draw[very thick] (120:3) -- (270:3);
    \draw \foreach \x in {0,30,...,330}
    {
      (\x:3)  node {}
    };
\end{tikzpicture}\quad
\\
\vspace{0.1in}
\begin{tikzpicture}[thick,scale=0.25]
    \draw (120:3) -- (240:3);
    \draw (120:3) -- (210:3);
    \draw (0:3) -- (30:3);
    \draw (210:3) -- (240:3);
    \draw (90:3) -- (120:3);
    \draw (90:3) -- (210:3);
    \draw (210:3) -- (330:3);
    \draw (330:3) -- (30:3);
    \draw (270:3) -- (300:3);
    \draw (120:3) -- (270:3);
    \draw[very thick] (30:3) -- (270:3);
   \draw \foreach \x in {0,30,...,330}
    {
      (\x:3)  node {}
    };
\end{tikzpicture}\quad
\begin{tikzpicture}[thick,scale=0.25]
    \draw (120:3) -- (240:3);
    \draw (120:3) -- (210:3);
    \draw (0:3) -- (30:3);
    \draw (210:3) -- (240:3);
    \draw (90:3) -- (120:3);
    \draw (90:3) -- (210:3);
    \draw (210:3) -- (330:3);
    \draw (330:3) -- (30:3);
    \draw (270:3) -- (300:3);
    \draw (120:3) -- (270:3);
    \draw (30:3) -- (270:3);
    \draw[very thick] (90:3) -- (300:3);
    \draw \foreach \x in {0,30,...,330}
    {
      (\x:3)  node {}
    };
\end{tikzpicture}\quad
\begin{tikzpicture}[thick,scale=0.25]
    \draw (120:3) -- (240:3);
    \draw (120:3) -- (210:3);
    \draw (0:3) -- (30:3);
    \draw (210:3) -- (240:3);
    \draw (90:3) -- (120:3);
    \draw (90:3) -- (210:3);
    \draw (210:3) -- (330:3);
    \draw (330:3) -- (30:3);
    \draw (270:3) -- (300:3);
    \draw (120:3) -- (270:3);
    \draw (30:3) -- (270:3);
    \draw (90:3) -- (300:3);
    \draw[very thick] (90:3) -- (30:3);
    \draw \foreach \x in {0,30,...,330}
    {
      (\x:3)  node {}
    };
\end{tikzpicture}\quad
\begin{tikzpicture}[thick,scale=0.25]
    \draw (120:3) -- (240:3);
    \draw (120:3) -- (210:3);
    \draw (0:3) -- (30:3);
    \draw (210:3) -- (240:3);
    \draw (90:3) -- (120:3);
    \draw (90:3) -- (210:3);
    \draw (210:3) -- (330:3);
    \draw (330:3) -- (30:3);
    \draw (270:3) -- (300:3);
    \draw (120:3) -- (270:3);
    \draw (30:3) -- (270:3);
    \draw (90:3) -- (300:3);
    \draw (90:3) -- (30:3);
    \draw[very thick] (30:3) -- (120:3);
    \draw \foreach \x in {0,30,...,330}
    {
      (\x:3)  node {}
    };
\end{tikzpicture}\quad
\begin{tikzpicture}[thick,scale=0.25]
    \draw (120:3) -- (240:3);
    \draw (120:3) -- (210:3);
    \draw (0:3) -- (30:3);
    \draw (210:3) -- (240:3);
    \draw (90:3) -- (120:3);
    \draw (90:3) -- (210:3);
    \draw (210:3) -- (330:3);
    \draw (330:3) -- (30:3);
    \draw (270:3) -- (300:3);
    \draw (120:3) -- (270:3);
    \draw (30:3) -- (270:3);
    \draw (90:3) -- (300:3);
    \draw (90:3) -- (30:3);
    \draw (30:3) -- (120:3);
    \draw[very thick] (90:3) -- (60:3);
    \draw \foreach \x in {0,30,...,330}
    {
      (\x:3)  node {}
    };
\end{tikzpicture}\quad
\\
\vspace{0.1in}
\begin{tikzpicture}[thick,scale=0.25]
    \draw (120:3) -- (240:3);
    \draw (120:3) -- (210:3);
    \draw (0:3) -- (30:3);
    \draw (210:3) -- (240:3);
    \draw (90:3) -- (120:3);
    \draw (90:3) -- (210:3);
    \draw (210:3) -- (330:3);
    \draw (330:3) -- (30:3);
    \draw (270:3) -- (300:3);
    \draw (120:3) -- (270:3);
    \draw (30:3) -- (270:3);
    \draw (90:3) -- (300:3);
    \draw (90:3) -- (30:3);
    \draw (30:3) -- (120:3);
    \draw (90:3) -- (60:3);
    \draw[very thick] (210:3) -- (300:3);
    \draw \foreach \x in {0,30,...,330}
    {
      (\x:3)  node {}
    };
\end{tikzpicture}\quad
\begin{tikzpicture}[thick,scale=0.25]
    \draw (120:3) -- (240:3);
    \draw (120:3) -- (210:3);
    \draw (0:3) -- (30:3);
    \draw (210:3) -- (240:3);
    \draw (90:3) -- (120:3);
    \draw (90:3) -- (210:3);
    \draw (210:3) -- (330:3);
    \draw (330:3) -- (30:3);
    \draw (270:3) -- (300:3);
    \draw (120:3) -- (270:3);
    \draw (30:3) -- (270:3);
    \draw (90:3) -- (300:3);
    \draw (90:3) -- (30:3);
    \draw (30:3) -- (120:3);
    \draw (90:3) -- (60:3);
    \draw (210:3) -- (300:3);
    \draw[very thick] (210:3) -- (60:3);
    \draw \foreach \x in {0,30,...,330}
    {
      (\x:3)  node {}
    };
\end{tikzpicture}\quad
\begin{tikzpicture}[thick,scale=0.25]
    \draw (120:3) -- (240:3);
    \draw (120:3) -- (210:3);
    \draw (0:3) -- (30:3);
    \draw (210:3) -- (240:3);
    \draw (90:3) -- (120:3);
    \draw (90:3) -- (210:3);
    \draw (210:3) -- (330:3);
    \draw (330:3) -- (30:3);
    \draw (270:3) -- (300:3);
    \draw (120:3) -- (270:3);
    \draw (30:3) -- (270:3);
    \draw (90:3) -- (300:3);
    \draw (90:3) -- (30:3);
    \draw (30:3) -- (120:3);
     \draw (90:3) -- (60:3);
    \draw (210:3) -- (300:3);
    \draw (210:3) -- (60:3);
    \draw[very thick] (240:3) -- (30:3);
    \draw \foreach \x in {0,30,...,330}
    {
      (\x:3)  node {}
    };
\end{tikzpicture}\quad
\begin{tikzpicture}[thick,scale=0.25]
    \draw (120:3) -- (240:3);
    \draw (120:3) -- (210:3);
    \draw (0:3) -- (30:3);
    \draw (210:3) -- (240:3);
    \draw (90:3) -- (120:3);
    \draw (90:3) -- (210:3);
    \draw (210:3) -- (330:3);
    \draw (330:3) -- (30:3);
    \draw (270:3) -- (300:3);
    \draw (120:3) -- (270:3);
    \draw (30:3) -- (270:3);
    \draw (90:3) -- (300:3);
    \draw (90:3) -- (30:3);
    \draw (30:3) -- (120:3);

    \draw (90:3) -- (60:3);
    \draw (210:3) -- (300:3);
    \draw (210:3) -- (60:3);

    \draw (240:3) -- (30:3);

    \draw[very thick] (120:3) -- (150:3);
    \draw \foreach \x in {0,30,...,330}
    {
      (\x:3)  node {}
    };
\end{tikzpicture}\quad
\begin{tikzpicture}[thick,scale=0.25]
    \draw (120:3) -- (240:3);
    \draw (120:3) -- (210:3);
    \draw (0:3) -- (30:3);
    \draw (210:3) -- (240:3);
    \draw (90:3) -- (120:3);
    \draw (90:3) -- (210:3);
    \draw (210:3) -- (330:3);
    \draw (330:3) -- (30:3);
    \draw (270:3) -- (300:3);
    \draw (120:3) -- (270:3);
    \draw (30:3) -- (270:3);
    \draw (90:3) -- (300:3);
    \draw (90:3) -- (30:3);
    \draw (30:3) -- (120:3);
    \draw (90:3) -- (60:3);
    \draw (210:3) -- (300:3);
    \draw (210:3) -- (60:3);
    \draw (240:3) -- (30:3);
    \draw (120:3) -- (150:3);
    \draw[very thick] (90:3) -- (180:3);
    \draw \foreach \x in {0,30,...,330}
    {
      (\x:3)  node {}
    };
\end{tikzpicture}\quad
\caption{The beginning of a random graph process on $n=12$ vertices.}
\label{fig:graph}
\end{figure}

It is also sometimes useful to consider a closely related ``random graph process''; see Figure \ref{fig:graph}. In the process $$\{ G(n,m) \}_{m=1}^{n \choose 2},$$
the $m$th edge is selected uniformly randomly from the remaining $${n \choose 2} - (m -1)$$ edges.

In random graph theory, one is usually interested in the asymptotic behavior of such graphs as $n \to \infty$, and $p=p(n)$.  We say that an event happens {\it with high probability} (abbreviated w.h.p.) if the probability approaches one as $n \to \infty$.
A celebrated theorem about the topology of random graphs is the following \cite{ER}.
 
\begin{theorem}[Erd\H{o}s--R\'enyi, 1959] \label{thm:ERS1} Let $\epsilon > 0$ be fixed, and $G \sim G(n,p)$.
If 
$$p \ge \frac{(1 + \epsilon) \log n }{ n},$$
then $G$ is connected w.h.p., and
if
$$p \le \frac{(1 - \epsilon) \log n }{ n},$$
then $G$ is disconnected w.h.p.
\end{theorem} 

(Although this is generally attributed to Erd\H{o}s and R\'enyi, it seems that Erd\H{o}s and R\'enyi actually proved the analogous statement for $G(n,m)$, and Stepanov first proved the $G(n,p)$ statement in 1969  \cite{S69}.)

The strongest possible statement of the Erd\H{o}s--R\'enyi theorem is slightly sharper than this, as follows.  Let $$\tilde{\beta_0}(G) = \# \mbox{ of connected components of } G - 1,$$ i.e.\ to a topologist, the reduced $0$th Betti number of $G$.

\begin{theorem}[Erd\H{o}s--R\'enyi, 1959] \label{thm:ERS2} Let $G = G(n,p)$ where $$ p = \frac{ \log n + c}{n},$$ and $c \in \R$ is fixed.
Then as $n \to \infty$, $\tilde{\beta}_0(G)$ is asymptotically Poisson distributed with mean $e^{-c}$.
In particular, 
$$ \prob [ G \mbox{ is connected} ]  \to e^{-e^{-c}}$$    
\end{theorem} 

\begin{corollary}  \label{cor:ER} Let $\omega \to \infty$ arbitrarily slowly as $n \to \infty$. 

If 
$$p \ge \frac{ \log n + \omega }{ n},$$
then $G$ is connected w.h.p., and
if
$$p \le \frac{\log n - \omega }{ n},$$
then $G$ is disconnected w.h.p.
\end{corollary}

Is there any way we could guess that $p = \log n / n$ the right threshold, if we did not already know that?  To give an intuition for this, we set $p = (\log n + c) / n$ with $c \in \R$ fixed, and ask what is the expected number of isolated vertices. The probability that a vertex is isolated is $(1-p)^{n-1}$, by independence. Then by linearity of expectation, the expected number of isolated vertices $X_0$ is given by
$$ \expect [ X_0]  = n (1-p)^{n-1}.$$
Since $p \to 0$ as $n \to \infty$, it is reasonable to approximate $1-p$ by $e^{-p}$, and then it follows easily that $\expect [ X_0]  \to e^{-c}$ as $n \to \infty$. With a little more work, for example by computing all the higher moments, one can show that in fact $X_0$ converges in law to a Poisson distribution with mean $e^{-c}$.  See, for example, Chapter 8 in \cite{Alon} for an overview of limit theorems in random graph theory.

To finish the proof of Theorem \ref{thm:ERS2}, one also needs a structure theorem, namely that for $p$ is in this range, w.h.p.\  $G \sim G(n,p)$ consists of only two kinds of connected components: a unique ``giant component,'' and isolated vertices.  Given this structure, $G$ is connected if and only if there are no isolated vertices. See Chapter 7 of \cite{Bollo} for a complete proof.

Corollary \ref{cor:ER} shows that $p = \log n / n$ is a {\it sharp threshold} for connectivity of $G(n,p)$, meaning that the phase transition from probability $0$ to probability $1$ happens within a very narrow window.  More precisely, a function $f$ is said to be a sharp threshold for a graph property $\mathcal{P}$ if there exists a function $g = o(f)$ such that 

\begin{displaymath}
\prob [ G(n,p) \in \mathcal{P} ]  \to \left\{
     \begin{array}{ll}
       1 & : p \ge f + g\\
       & \\
       0 & : p \le f-g 
     \end{array}
   \right.
\end{displaymath} 

Exactly which monotone graph properties have sharp thresholds is a question which has been extensively studied. See for example the paper of Friedghut with appendix by Bourgain \cite{FB99}.

Corollary \ref{cor:ER} can be summarized as saying that the threshold function for ``$G$ is connected'' is the same as the threshold function for ``$G$ has no isolated vertices.''  The following result of Bollob\'as and Thomasson takes this idea all the way to its logical conclusion \cite{bt}.

\begin{theorem}\label{thm:BT} For a random graph process $\{G(n,m)\}_{m=1}^{n \choose 2}$, with high probability
$$
\min\{M: G(n,M) \text{ has no isolated vertices} \}=\min\{M:G(n,M)\text{ is connected}\}.
$$
\label{thm:BT}
\end{theorem}

We leave it to readers to convince themselves that Theorem \ref{thm:BT} is even sharper than Corollary \ref{cor:ER}.

\medskip

There is another topological phase transition for $G \sim G(n,p)$, namely where cycles first appear, or to a topologist, where $H_1 (G)$ first becomes nontrivial.  See Pittel \cite{P88} for a proof of the following characterization of the appearance of cycles.

\begin{theorem} \label{thm:cycle}  Let $G \sim G(n,p)$, where $p = c/ n$ and $c > 0$ is constant.
\begin{displaymath}
\prob [ H_1(G) \neq 0 ]  \to \left\{
     \begin{array}{lr}
       1 & : c \ge 1\\
       & \\
 \sqrt{1-c} \, \, \exp(c/2 + c^2/4) & : c < 1
     \end{array}
   \right.
\end{displaymath} 
\end{theorem}

Note that in contrast to the connectivity threshold, this threshold is sharp on one side but not on the other.

\section{Random $2$-complexes} \label{sec:2com}

\begin{figure}
\tikzstyle{every node}=[circle, draw, fill=orange,
                        inner sep=0pt, minimum width=4pt]
\begin{tikzpicture}[scale=0.25]
\draw \foreach \a in {0,30,...,330}
{ \foreach \b in {0, 30, ..., 330}
{
(\a:3) -- (\b:3)
}
};

    \draw \foreach \x in {0,30,...,330}
    {
      (\x:3)  node {}
    };
\end{tikzpicture}\quad
\begin{tikzpicture}[scale=0.25]
\draw [black, fill=blue, opacity = 0.3] (90:3) -- (240:3) -- (270:3) -- cycle;
\draw \foreach \a in {0,30,...,330}
{ \foreach \b in {0, 30, ..., 330}
{
(\a:3) -- (\b:3)
}
};
\draw[very thick] (90:3)--(240:3)--(270:3)--cycle;
    \draw \foreach \x in {0,30,...,330}
    {
      (\x:3)  node {}
    };
\end{tikzpicture}\quad
\begin{tikzpicture}[scale=0.25]
\draw [black, fill=blue, opacity = 0.3] (90:3) -- (240:3) -- (270:3) -- cycle;
\draw [black, fill=blue, opacity = 0.3] (0:3) -- (120:3) -- (330:3) -- cycle;
\draw \foreach \a in {0,30,...,330}
{ \foreach \b in {0, 30, ..., 330}
{
(\a:3) -- (\b:3)
}
};
\draw[very thick] (0:3)--(120:3)--(330:3)--cycle;
    \draw \foreach \x in {0,30,...,330}
    {
      (\x:3)  node {}
    };
\end{tikzpicture}\quad
\begin{tikzpicture}[scale=0.25]
\draw [black, fill=blue, opacity = 0.3] (90:3) -- (240:3) -- (270:3) -- cycle;
\draw [black, fill=blue, opacity = 0.3] (0:3) -- (120:3) -- (330:3) -- cycle;
\draw [black, fill=blue, opacity = 0.3] (150:3) -- (120:3) -- (300:3) -- cycle;
\draw \foreach \a in {0,30,...,330}
{ \foreach \b in {0, 30, ..., 330}
{
(\a:3) -- (\b:3)
}
};
\draw[very thick] (150:3)--(120:3)--(300:3)--cycle;
    \draw \foreach \x in {0,30,...,330}
    {
      (\x:3)  node {}
    };
\end{tikzpicture}\quad
\begin{tikzpicture}[scale=0.25]
\draw [black, fill=blue, opacity = 0.3] (90:3) -- (240:3) -- (270:3) -- cycle;
\draw [black, fill=blue, opacity = 0.3] (0:3) -- (120:3) -- (330:3) -- cycle;
\draw [black, fill=blue, opacity = 0.3] (150:3) -- (120:3) -- (300:3) -- cycle;
\draw [black, fill=blue, opacity = 0.3] (120:3) -- (210:3) -- (330:3) -- cycle;
\draw \foreach \a in {0,30,...,330}
{ \foreach \b in {0, 30, ..., 330}
{
(\a:3) -- (\b:3)
}
};
\draw[very thick] (120:3)--(210:3)--(330:3)--cycle;
    \draw \foreach \x in {0,30,...,330}
    {
      (\x:3)  node {}
    };
\end{tikzpicture}\quad
\\
\vspace{0.1in}
\begin{tikzpicture}[scale=0.25]
\draw [black, fill=blue, opacity = 0.3] (90:3) -- (240:3) -- (270:3) -- cycle;
\draw [black, fill=blue, opacity = 0.3] (0:3) -- (120:3) -- (330:3) -- cycle;
\draw [black, fill=blue, opacity = 0.3] (150:3) -- (120:3) -- (300:3) -- cycle;
\draw [black, fill=blue, opacity = 0.3] (120:3) -- (210:3) -- (330:3) -- cycle;
\draw [black, fill=blue, opacity = 0.3] (30:3) -- (120:3) -- (150:3) -- cycle;
\draw \foreach \a in {0,30,...,330}
{ \foreach \b in {0, 30, ..., 330}
{
(\a:3) -- (\b:3)
}
};
\draw[very thick] (30:3)--(120:3)--(150:3)--cycle;
    \draw \foreach \x in {0,30,...,330}
    {
      (\x:3)  node {}
    };
\end{tikzpicture}\quad
\begin{tikzpicture}[scale=0.25]
\draw [black, fill=blue, opacity = 0.3] (90:3) -- (240:3) -- (270:3) -- cycle;
\draw [black, fill=blue, opacity = 0.3] (0:3) -- (120:3) -- (330:3) -- cycle;
\draw [black, fill=blue, opacity = 0.3] (150:3) -- (120:3) -- (300:3) -- cycle;
\draw [black, fill=blue, opacity = 0.3] (120:3) -- (210:3) -- (330:3) -- cycle;
\draw [black, fill=blue, opacity = 0.3] (30:3) -- (120:3) -- (150:3) -- cycle;
\draw [black, fill=blue, opacity = 0.3] (60:3) -- (240:3) -- (330:3) -- cycle;
\draw \foreach \a in {0,30,...,330}
{ \foreach \b in {0, 30, ..., 330}
{
(\a:3) -- (\b:3)
}
};
\draw[very thick] (60:3)--(240:3)--(330:3)--cycle;
    \draw \foreach \x in {0,30,...,330}
    {
      (\x:3)  node {}
    };
\end{tikzpicture}\quad
\begin{tikzpicture}[scale=0.25]
\draw [black, fill=blue, opacity = 0.3] (90:3) -- (240:3) -- (270:3) -- cycle;
\draw [black, fill=blue, opacity = 0.3] (0:3) -- (120:3) -- (330:3) -- cycle;
\draw [black, fill=blue, opacity = 0.3] (150:3) -- (120:3) -- (300:3) -- cycle;
\draw [black, fill=blue, opacity = 0.3] (120:3) -- (210:3) -- (330:3) -- cycle;
\draw [black, fill=blue, opacity = 0.3] (30:3) -- (120:3) -- (150:3) -- cycle;
\draw [black, fill=blue, opacity = 0.3] (60:3) -- (240:3) -- (330:3) -- cycle;
\draw [black, fill=blue, opacity = 0.3] (30:3) -- (90:3) -- (240:3) -- cycle;

\draw \foreach \a in {0,30,...,330}
{ \foreach \b in {0, 30, ..., 330}
{
(\a:3) -- (\b:3)
}
};
\draw[very thick] (30:3)--(90:3)--(240:3)--cycle;
    \draw \foreach \x in {0,30,...,330}
    {
      (\x:3)  node {}
    };
\end{tikzpicture}\quad
\begin{tikzpicture}[scale=0.25]
\draw [black, fill=blue, opacity = 0.3] (90:3) -- (240:3) -- (270:3) -- cycle;
\draw [black, fill=blue, opacity = 0.3] (0:3) -- (120:3) -- (330:3) -- cycle;
\draw [black, fill=blue, opacity = 0.3] (150:3) -- (120:3) -- (300:3) -- cycle;
\draw [black, fill=blue, opacity = 0.3] (120:3) -- (210:3) -- (330:3) -- cycle;
\draw [black, fill=blue, opacity = 0.3] (30:3) -- (120:3) -- (150:3) -- cycle;
\draw [black, fill=blue, opacity = 0.3] (60:3) -- (240:3) -- (330:3) -- cycle;
\draw [black, fill=blue, opacity = 0.3] (30:3) -- (90:3) -- (240:3) -- cycle;
\draw [black, fill=blue, opacity = 0.3] (60:3) -- (90:3) -- (240:3) -- cycle;

\draw \foreach \a in {0,30,...,330}
{ \foreach \b in {0, 30, ..., 330}
{
(\a:3) -- (\b:3)
}
};
\draw[very thick] (60:3)--(90:3)--(240:3)--cycle;
    \draw \foreach \x in {0,30,...,330}
    {
      (\x:3)  node {}
    };
\end{tikzpicture}\quad
\begin{tikzpicture}[scale=0.25]
\draw [black, fill=blue, opacity = 0.3] (90:3) -- (240:3) -- (270:3) -- cycle;
\draw [black, fill=blue, opacity = 0.3] (0:3) -- (120:3) -- (330:3) -- cycle;
\draw [black, fill=blue, opacity = 0.3] (150:3) -- (120:3) -- (300:3) -- cycle;
\draw [black, fill=blue, opacity = 0.3] (120:3) -- (210:3) -- (330:3) -- cycle;
\draw [black, fill=blue, opacity = 0.3] (30:3) -- (120:3) -- (150:3) -- cycle;
\draw [black, fill=blue, opacity = 0.3] (60:3) -- (240:3) -- (330:3) -- cycle;
\draw [black, fill=blue, opacity = 0.3] (30:3) -- (90:3) -- (240:3) -- cycle;
\draw [black, fill=blue, opacity = 0.3] (60:3) -- (90:3) -- (240:3) -- cycle;
\draw [black, fill=blue, opacity = 0.3] (120:3) -- (90:3) -- (330:3) -- cycle;

\draw \foreach \a in {0,30,...,330}
{ \foreach \b in {0, 30, ..., 330}
{
(\a:3) -- (\b:3)
}
};
\draw[very thick] (120:3)--(90:3)--(330:3)--cycle;
    \draw \foreach \x in {0,30,...,330}
    {
      (\x:3)  node {}
    };
\end{tikzpicture}\quad
\\
\vspace{0.1in}
\begin{tikzpicture}[scale=0.25]
\draw [black, fill=blue, opacity = 0.3] (90:3) -- (240:3) -- (270:3) -- cycle;
\draw [black, fill=blue, opacity = 0.3] (0:3) -- (120:3) -- (330:3) -- cycle;
\draw [black, fill=blue, opacity = 0.3] (150:3) -- (120:3) -- (300:3) -- cycle;
\draw [black, fill=blue, opacity = 0.3] (120:3) -- (210:3) -- (330:3) -- cycle;
\draw [black, fill=blue, opacity = 0.3] (30:3) -- (120:3) -- (150:3) -- cycle;
\draw [black, fill=blue, opacity = 0.3] (60:3) -- (240:3) -- (330:3) -- cycle;
\draw [black, fill=blue, opacity = 0.3] (30:3) -- (90:3) -- (240:3) -- cycle;
\draw [black, fill=blue, opacity = 0.3] (60:3) -- (90:3) -- (240:3) -- cycle;
\draw [black, fill=blue, opacity = 0.3] (120:3) -- (90:3) -- (330:3) -- cycle;
\draw [black, fill=blue, opacity = 0.3] (0:3) -- (90:3) -- (150:3) -- cycle;
\draw \foreach \a in {0,30,...,330}
{ \foreach \b in {0, 30, ..., 330}
{
(\a:3) -- (\b:3)
}
};
\draw[very thick] (0:3)--(90:3)--(150:3)--cycle;
    \draw \foreach \x in {0,30,...,330}
    {
      (\x:3)  node {}
    };
\end{tikzpicture}\quad
\begin{tikzpicture}[scale=0.25]
\draw [black, fill=blue, opacity = 0.3] (90:3) -- (240:3) -- (270:3) -- cycle;
\draw [black, fill=blue, opacity = 0.3] (0:3) -- (120:3) -- (330:3) -- cycle;
\draw [black, fill=blue, opacity = 0.3] (150:3) -- (120:3) -- (300:3) -- cycle;
\draw [black, fill=blue, opacity = 0.3] (120:3) -- (210:3) -- (330:3) -- cycle;
\draw [black, fill=blue, opacity = 0.3] (30:3) -- (120:3) -- (150:3) -- cycle;
\draw [black, fill=blue, opacity = 0.3] (60:3) -- (240:3) -- (330:3) -- cycle;
\draw [black, fill=blue, opacity = 0.3] (30:3) -- (90:3) -- (240:3) -- cycle;
\draw [black, fill=blue, opacity = 0.3] (60:3) -- (90:3) -- (240:3) -- cycle;
\draw [black, fill=blue, opacity = 0.3] (120:3) -- (90:3) -- (330:3) -- cycle;
\draw [black, fill=blue, opacity = 0.3] (0:3) -- (90:3) -- (150:3) -- cycle;
\draw [black, fill=blue, opacity = 0.3] (150:3) -- (210:3) -- (330:3) -- cycle;
\draw \foreach \a in {0,30,...,330}
{ \foreach \b in {0, 30, ..., 330}
{
(\a:3) -- (\b:3)
}
};
\draw[very thick] (150:3)--(210:3)--(330:3)--cycle;
    \draw \foreach \x in {0,30,...,330}
    {
      (\x:3)  node {}
    };
\end{tikzpicture}\quad
\begin{tikzpicture}[scale=0.25]
\draw [black, fill=blue, opacity = 0.3] (90:3) -- (240:3) -- (270:3) -- cycle;
\draw [black, fill=blue, opacity = 0.3] (0:3) -- (120:3) -- (330:3) -- cycle;
\draw [black, fill=blue, opacity = 0.3] (150:3) -- (120:3) -- (300:3) -- cycle;
\draw [black, fill=blue, opacity = 0.3] (120:3) -- (210:3) -- (330:3) -- cycle;
\draw [black, fill=blue, opacity = 0.3] (30:3) -- (120:3) -- (150:3) -- cycle;
\draw [black, fill=blue, opacity = 0.3] (60:3) -- (240:3) -- (330:3) -- cycle;
\draw [black, fill=blue, opacity = 0.3] (30:3) -- (90:3) -- (240:3) -- cycle;
\draw [black, fill=blue, opacity = 0.3] (60:3) -- (90:3) -- (240:3) -- cycle;
\draw [black, fill=blue, opacity = 0.3] (120:3) -- (90:3) -- (330:3) -- cycle;
\draw [black, fill=blue, opacity = 0.3] (0:3) -- (90:3) -- (150:3) -- cycle;
\draw [black, fill=blue, opacity = 0.3] (150:3) -- (210:3) -- (330:3) -- cycle;
\draw [black, fill=blue, opacity = 0.3] (120:3) -- (240:3) -- (300:3) -- cycle;
\draw \foreach \a in {0,30,...,330}
{ \foreach \b in {0, 30, ..., 330}
{
(\a:3) -- (\b:3)
}
};
\draw[very thick] (120:3)--(240:3)--(300:3)--cycle;
    \draw \foreach \x in {0,30,...,330}
    {
      (\x:3)  node {}
    };
\end{tikzpicture}\quad
\begin{tikzpicture}[scale=0.25]
\draw [black, fill=blue, opacity = 0.3] (90:3) -- (240:3) -- (270:3) -- cycle;
\draw [black, fill=blue, opacity = 0.3] (0:3) -- (120:3) -- (330:3) -- cycle;
\draw [black, fill=blue, opacity = 0.3] (150:3) -- (120:3) -- (300:3) -- cycle;
\draw [black, fill=blue, opacity = 0.3] (120:3) -- (210:3) -- (330:3) -- cycle;
\draw [black, fill=blue, opacity = 0.3] (30:3) -- (120:3) -- (150:3) -- cycle;
\draw [black, fill=blue, opacity = 0.3] (60:3) -- (240:3) -- (330:3) -- cycle;
\draw [black, fill=blue, opacity = 0.3] (30:3) -- (90:3) -- (240:3) -- cycle;
\draw [black, fill=blue, opacity = 0.3] (60:3) -- (90:3) -- (240:3) -- cycle;
\draw [black, fill=blue, opacity = 0.3] (120:3) -- (90:3) -- (330:3) -- cycle;
\draw [black, fill=blue, opacity = 0.3] (0:3) -- (90:3) -- (150:3) -- cycle;
\draw [black, fill=blue, opacity = 0.3] (150:3) -- (210:3) -- (330:3) -- cycle;
\draw [black, fill=blue, opacity = 0.3] (120:3) -- (240:3) -- (300:3) -- cycle;
\draw [black, fill=blue, opacity = 0.3] (30:3) -- (90:3) -- (150:3) -- cycle;
\draw \foreach \a in {0,30,...,330}
{ \foreach \b in {0, 30, ..., 330}
{
(\a:3) -- (\b:3)
}
};
\draw[very thick] (30:3)--(90:3)--(150:3)--cycle;
    \draw \foreach \x in {0,30,...,330}
    {
      (\x:3)  node {}
    };
\end{tikzpicture}\quad
\begin{tikzpicture}[scale=0.25]
\draw [black, fill=blue, opacity = 0.3] (90:3) -- (240:3) -- (270:3) -- cycle;
\draw [black, fill=blue, opacity = 0.3] (0:3) -- (120:3) -- (330:3) -- cycle;
\draw [black, fill=blue, opacity = 0.3] (150:3) -- (120:3) -- (300:3) -- cycle;
\draw [black, fill=blue, opacity = 0.3] (120:3) -- (210:3) -- (330:3) -- cycle;
\draw [black, fill=blue, opacity = 0.3] (30:3) -- (120:3) -- (150:3) -- cycle;
\draw [black, fill=blue, opacity = 0.3] (60:3) -- (240:3) -- (330:3) -- cycle;
\draw [black, fill=blue, opacity = 0.3] (30:3) -- (90:3) -- (240:3) -- cycle;
\draw [black, fill=blue, opacity = 0.3] (60:3) -- (90:3) -- (240:3) -- cycle;
\draw [black, fill=blue, opacity = 0.3] (120:3) -- (90:3) -- (330:3) -- cycle;
\draw [black, fill=blue, opacity = 0.3] (0:3) -- (90:3) -- (150:3) -- cycle;
\draw [black, fill=blue, opacity = 0.3] (150:3) -- (210:3) -- (330:3) -- cycle;
\draw [black, fill=blue, opacity = 0.3] (120:3) -- (240:3) -- (300:3) -- cycle;
\draw [black, fill=blue, opacity = 0.3] (30:3) -- (90:3) -- (150:3) -- cycle;
\draw [black, fill=blue, opacity = 0.3] (60:3) -- (300:3) -- (330:3) -- cycle;
\draw \foreach \a in {0,30,...,330}
{ \foreach \b in {0, 30, ..., 330}
{
(\a:3) -- (\b:3)
}
};
\draw[very thick] (60:3)--(300:3)--(330:3)--cycle;
    \draw \foreach \x in {0,30,...,330}
    {
      (\x:3)  node {}
    };
\end{tikzpicture}\quad
\caption{The beginning of a random $2$-complex process on $n=12$ vertices.}
\label{fig:2com}
\end{figure}
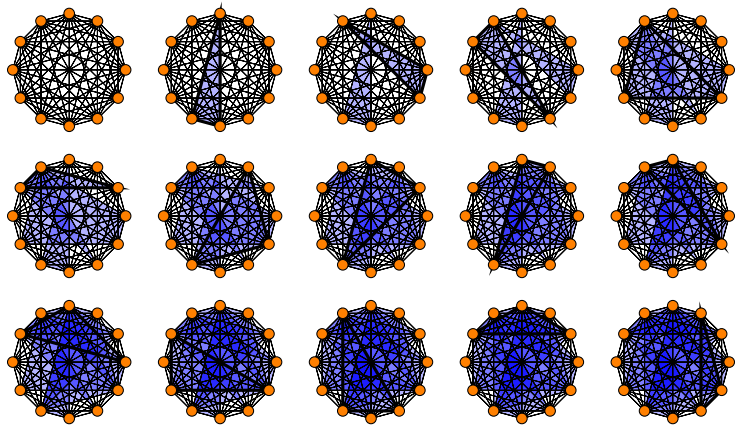

Linial and Meshulam initiated the topological study of random $2$-dimensional simplicial complexes $Y(n,p)$ in \cite{LM}.  This model random simplicial complex is defined to have vertex set $[n]$, edge set $[n] \choose 2$ (i.e.\ the underlying graph is a complete graph), and each of the $n \choose 3$ possible triangle faces is included with probability $p$, independently.  The beginning of a random $2$-complex process is illustrated in Figure \ref{fig:2com}.

The main result of \cite{LM} is a cohomological analogue of Theorem \ref{thm:ERS1}.

\begin{theorem}[Linial--Meshulam, 2006] \label{thm:LM} Let $\epsilon > 0$ be fixed and $Y \sim Y(n,p)$. If
$$p \ge (2 + \epsilon) \log n / n,$$
then w.h.p.\ $H^1(Y, \Z / 2) = 0$, and if
$$p \le (2 - \epsilon) \log n / n,$$
then w.h.p.\ $H^1(Y, \Z / 2) \neq 0$.
\end{theorem} 

Although Theorem \ref{thm:LM} is analogous to Theorem \ref{thm:ERS1}, the proof is much harder. Linial and Meshulam used a combination of co-isoperimetric inequalities (analogous to Cheeger constant of a graph) and intricate cocycle-counting arguments.

A few comments might be in order.

\begin{enumerate}
\item The Linial--Meshulam theorem is sharper than this, analogous to Corollary \ref{cor:ER}, but in this survey article we sometimes trade the strongest result for a simpler statement.
\item The threshold $p = 2 \log n / n$ is just what is necessary in order to ensure that there are no isolated edges. The analogue of the ``stopping time'' Theorem \ref{thm:BT} for the random $2$-complex process was recently established in \cite{Pittel}. 
\item Theorem \ref{thm:LM} is stated for cohomology, but universal coefficients for homology and cohomology give the analogous result for homology.
\item  It was shown by Meshulam and Wallach \cite{MW} that the same result holds with $(\Z / \ell)$-coefficients for every fixed $\ell$ (and for an analogous $d$-dimensional model).
\item The threshold for homology with $\Z$-coefficients is still unknown.  It might seem that this would follow from Meshulam and Wallach's work, but the problem is that there could be $\ell$-torsion, where $\ell$ is growing with $n$.
\end{enumerate}

Linial and Meshulam asked in \cite{LM} for the threshold for simple connectivity, and this was eventually shown to be much larger \cite{BHK}.

\begin{theorem}[Babson et al., 2011] \label{thm:BHK} Let $\epsilon > 0$ be fixed and $Y \sim Y(n,p)$.  If
$$p \ge \frac{n^{ \epsilon}}{\sqrt{n}},$$
then w.h.p.\ $\pi_1(Y) = 0$, and if
$$p \le \frac{n^{-\epsilon}}{\sqrt{n}},$$
then w.h.p.\ $\pi_1(Y)$ is a nontrivial group, hyperbolic in the sense of Gromov.
\end{theorem} 

The study of fundamental groups of random $2$-complexes is continued in \cite{hkp}.  A group $G$ is said to have Kazhdan's property (T) if the trivial representation is an isolated point in the unitary dual of $G$ equipped with the Fell topology.  More intuitively, Property (T) is an ``expander-like'' property of groups; in fact the first explicit examples of expanders were constructed by Margulis from Cayley graphs of quotients of (T) groups such as $SL(3,\Z)$ \cite{Margulis73}. For a comprehensive overview of Property (T) see the monograph \cite{BdV}.

The following theorem shows that the threshold for $\pi_1(Y)$ to have Kazhdan's Property (T) is the same as the vanishing threshold for $H^1(Y,\Z/2)$ \cite{hkp}.

\begin{theorem}[Hoffman et al., 2012] \label{thm:HKP} Let $\epsilon > 0$ be fixed and $Y \sim Y(n,p)$.  Then as $n \to \infty$,
\begin{displaymath}
\prob [ \pi_1( Y) \mbox{ is Kazhdan} ]  \to \left\{
     \begin{array}{lr}
       1 & : p \ge (2 + \epsilon) \log n / n\\
       & \\
       0 & : p \le (2 - \epsilon) \log n / n
     \end{array}
   \right.
\end{displaymath} 
\end{theorem}

The proof that $Y(n,p)$ is Kazhdan when $p \ge (2 + \epsilon) \log n / n$ utilizes the following theorem of \.Zuk \cite{Zuk}.

\begin{theorem}[\.Zuk] \label{thm:tool}
If $X$ is a pure $2$-dimensional locally-finite simplicial complex so that for every vertex $v$, the vertex link $\link_v(X)$ is connected and the normalized graph Laplacian $L= L(\link_v(X))$ has smallest positive eigenvalue $\lambda_2(L) > 1/2$, then $\pi_1(X)$ has property (T).
\end{theorem}

The link of a vertex in the random $2$-complex $Y(n,p)$ has the same probability distribution as a random graph $G(n-1,p)$.  So \.Zuk's theorem reduces the proof of Theorem \ref{thm:HKP} to a question about Laplacians of random graphs.  However, new results about such Laplacians are still required in order to prove Theorem \ref{thm:HKP}.  Establishing the following comprises most of the work in \cite{hkp}.
 
\begin{theorem}[Hoffman et al., 2012]
\label{thm:ergap1}
Fix $k \ge 0$ and $\epsilon >0$.  Let $0=\lambda_1 \leq \lambda_2 \leq \cdots \leq \lambda_n \le 2$ be the eigenvalues of the normalized Laplacian of the random graph $G(n,p)$.  There is a constant $C=C(k)$ so that when
$$p \ge  \frac{ (k+1) \log n + C \sqrt{ \log n } \log\log n}{n}$$
is satisfied, then
$$\lambda_2 > 1-\epsilon,$$ with probability at least $1 - o(n^{-k})$.
\end{theorem}

The proof of Theorem \ref{thm:HKP} only requires Theorem \ref{thm:ergap1} with $k=1$ and $\epsilon = 1/2$.  By \.Zuk's theorem, the proof is almost immediate once Theorem \ref{thm:ergap1} is established:  There are $n$ vertex links, and for each one, the probability that its spectral gap is smaller than $1/2$ is $o(n^{-1})$.  So the probability that there is at least one vertex link with spectral gap too small is $o(1)$ by a union bound  --- the probability that at least one bad event occurs is never more than the sum of the probabilities of the individual bad events.

The $k=0$ case of Theorem \ref{thm:ergap1} is also of particular interest, as this gets very close to the connectivity threshold for $G(n,p)$.  It would be interesting to see just how close one can get.  Consider a random graph process, for example, adding one edge at a time: is the graph already an expander (smallest positive eigenvalue bounded away from zero) at the moment of connectivity? We discuss applications of Theorem \ref{thm:ergap1} with other values of $k$ in Section \ref{sec:rfc}.
 
\bigskip 
 
Both Theorem \ref{thm:LM} and Theorem \ref{thm:HKP} have the following corollary.

\begin{corollary} \label{thm:Q2} Let $\epsilon > 0$ be fixed and $Y \sim Y(n,p)$.  If
$$p \ge (2 + \epsilon) \log n / n,$$
then w.h.p.\ $H_1(Y, \Q) = 0$,
and if 
$$p \le (2 - \epsilon) \log n / n,$$
then w.h.p.\ $H_1(Y, \Q \neq 0$.
\end{corollary} 
 
This follows from Theorem \ref{thm:LM} by the universal coefficient theorem, and follows from Theorem \ref{thm:HKP}, since a Kazhdan group has finite abelianization.

\bigskip

Costa and Farber described two phase transitions for the cohomological dimension $\cd \pi_1(Y)$ \cite{CF12,CF13}.
\begin{theorem} [Costa--Farber] Let $Y \sim Y(n,p)$, and set $p = n^{-\alpha}$.
\begin{enumerate}
\item If $\alpha > 1$ then w.h.p.\ $\cd \pi_1(Y) = 1$ ,
\item if $3/5 < \alpha < 1$ then w.h.p.\ $\cd  \pi_1(Y)=2$, and
\item if $1/2 < \alpha < 3/5$ then w.h.p.\ $\cd \pi_1(Y) = \infty$.
\end{enumerate}
\end{theorem}

To give some intuition for where the exponent $3/5$ comes from: the smallest triangulation of the projective plane has $6$ vertices and $10$ faces, and these appear as subcomplexes once $p \gg n^{-6/10} = n^{-3/5}$. Costa and Farber show that the induced map on fundamental group is injective, hence for $p$ in this regime there is $2$-torsion in $\pi_1(Y)$. 

\bigskip

The random fundamental groups are related to other models random finitely presented groups, as in the survey article \cite{Ollivier}, especially the ``triangular'' model studied earlier by \.Zuk, but the geometric and topological considerations in the random fundamental group setting are somewhat more complicated.

As an aside, very different models of random groups based on random graphs, with large but finite cohomological dimension, for example, have recently been studied by Charney--Farber \cite{RRAAG}, and by Davis--Kahle \cite{dual}.

\bigskip

The two-dimensional analogue of Theorem \ref{thm:cycle} describing the first appearance of cycles has also been studied, in a series of papers by Kozlov \cite{Kozlov}, Cohen et al.\ \cite{CCFK12}, and Aronshtam--Linial \cite{AL}. In this last paper, the following strong estimate is established.

\begin{theorem}[Aronshtam--Linial, 2012]
\label{thm:2cycle}  Let $Y \sim Y(n,p)$ where $p = c/ n$ and $c > 0$ is constant.
If $c > 2.75381\dots$ then w.h.p. 
$H_2(Y) \neq 0.$
\end{theorem}

Cohen et al.\ had already noticed that this theorem holds with $c > 3$, essentially for linear algebraic reasons. It is surprisingly difficult to improve the constant past $3$. The number $2.75381\dots$ is the unique positive solution of a fairly complicated equation described in the paper, and Aronshtam and Linial conjecture that this constant is best possible.

\bigskip

Homology vanishing,  simple connectivity, and Property (T) are all monotone properties for random $2$-complexes. This is in contrast to what we will see in Section \ref{sec:rfc}, where each homology group $H_k(X)$,$k \ge 1$, passes through two distinct phase transitions: one where it first appears, and a later one where it vanishes.
 
\section{Random flag complexes} \label{sec:rfc}

The {\it flag complex} $X(H)$ of a graph $H$ is the maximal simplicial complex compatible with $H$ as its $1$-skeleton; in other words, the $i$-dimensional faces of $X(H)$ correspond to the cliques of order $i+1$ in $H$. (Such complexes have apparently arisen independently several times, and $X(H)$ is also sometimes called the {\it clique complex} or the {\it Vietoris--Rips complex} of $H$.)

We define the random flag complex $X(n,p)$ to be the flag complex of the random graph $G(n,p)$.  Every simplicial complex is homeomorphic to a flag complex, e.g.\ by taking the barycentric subdivision.  So $X(n,p)$ puts a measure on a wide variety of topologies as $n \to \infty$.

\begin{figure}
\tikzstyle{every node}=[circle, draw, fill=orange,
                        inner sep=0pt, minimum width=4pt]
\begin{tikzpicture}[scale=0.25]
    \draw \foreach \x in {0,30,...,330}
    {
      (\x:3)  node {}
    };
\end{tikzpicture}\quad
\begin{tikzpicture}[thick,scale=0.25]
    \draw[very thick](120:3) -- (240:3);
    \draw \foreach \x in {0,30,...,330}
    {
      (\x:3)  node {}
    };
\end{tikzpicture}\quad
\begin{tikzpicture}[thick,scale=0.25]
    \draw (120:3) -- (240:3);
    \draw[very thick] (120:3) -- (210:3);
    \draw \foreach \x in {0,30,...,330}
    {
      (\x:3)  node {}
    };
\end{tikzpicture}\quad
\begin{tikzpicture}[thick,scale=0.25]
    \draw (120:3) -- (240:3);
    \draw (120:3) -- (210:3);
    \draw[very thick] (0:3) -- (30:3);
    \draw \foreach \x in {0,30,...,330}
    {
      (\x:3)  node {}
    };
\end{tikzpicture}\quad
\begin{tikzpicture}[thick,scale=0.25]
    \draw (120:3) -- (240:3);
    \draw (120:3) -- (210:3);
    \draw (0:3) -- (30:3);
    \draw[very thick] (210:3) -- (240:3);
    \draw [black, fill=blue, opacity = 0.35] (210:3) -- (240:3) -- (120:3) -- cycle;
    \draw \foreach \x in {0,30,...,330}
    {
      (\x:3)  node {}
    };
\end{tikzpicture}\quad
\\
\vspace{0.1in}

\begin{tikzpicture}[thick,scale=0.25]
    \draw (120:3) -- (240:3);
    \draw (120:3) -- (210:3);
    \draw (0:3) -- (30:3);
    \draw  (210:3) -- (240:3);
    \draw [black, fill=blue, opacity = 0.35] (210:3) -- (240:3) -- (120:3) -- cycle;
    \draw[very thick] (90:3) -- (120:3);
    \draw \foreach \x in {0,30,...,330}
    {
      (\x:3)  node {}
    };
\end{tikzpicture}\quad
\begin{tikzpicture}[thick,scale=0.25]
    \draw (120:3) -- (240:3);
    \draw (120:3) -- (210:3);
    \draw (0:3) -- (30:3);
    \draw (210:3) -- (240:3);
    \draw [black, fill=blue, opacity=0.25] (210:3) -- (240:3) -- (120:3) -- cycle;
    \draw (90:3) -- (120:3);
    \draw[very thick]  (90:3) -- (210:3);
    \draw [black, fill=blue, opacity=0.25] (210:3) -- (90:3) -- (120:3) -- cycle;
    \draw \foreach \x in {0,30,...,330}
    {
      (\x:3)  node {}
    };
\end{tikzpicture}\quad
\begin{tikzpicture}[thick,scale=0.25]
    \draw (120:3) -- (240:3);
    \draw (120:3) -- (210:3);
    \draw (0:3) -- (30:3);
    \draw (210:3) -- (240:3);
   \draw [black, fill=blue, opacity=0.25] (210:3) -- (240:3) -- (120:3) -- cycle;
    \draw (90:3) -- (120:3);
    \draw (90:3) -- (210:3);
    \draw [black, fill=blue, opacity=0.25] (210:3) -- (90:3) -- (120:3) -- cycle;
    \draw[very thick] (210:3) -- (330:3);
    \draw \foreach \x in {0,30,...,330}
    {
      (\x:3)  node {}
    };
\end{tikzpicture}\quad
\begin{tikzpicture}[thick,scale=0.25]
    \draw (120:3) -- (240:3);
    \draw (120:3) -- (210:3);
    \draw (0:3) -- (30:3);
    \draw (210:3) -- (240:3);
  \draw [black, fill=blue, opacity=0.25] (210:3) -- (240:3) -- (120:3) -- cycle;
    \draw (90:3) -- (120:3);
    \draw (90:3) -- (210:3);
    \draw [black, fill=blue, opacity=0.25] (210:3) -- (90:3) -- (120:3) -- cycle;
    \draw (210:3) -- (330:3);
    \draw[very thick] (330:3) -- (30:3);
       \draw \foreach \x in {0,30,...,330}
    {
      (\x:3)  node {}
    };
\end{tikzpicture}\quad
\begin{tikzpicture}[thick,scale=0.25]
    \draw (120:3) -- (240:3);
    \draw (120:3) -- (210:3);
    \draw (0:3) -- (30:3);
    \draw (210:3) -- (240:3);
  \draw [black, fill=blue, opacity=0.25] (210:3) -- (240:3) -- (120:3) -- cycle;
    \draw (90:3) -- (120:3);
    \draw (90:3) -- (210:3);
    \draw [black, fill=blue, opacity=0.25] (210:3) -- (90:3) -- (120:3) -- cycle;
    \draw (210:3) -- (330:3);
    \draw (330:3) -- (30:3);
    \draw[very thick] (270:3) -- (300:3);
    \draw \foreach \x in {0,30,...,330}
    {
      (\x:3)  node {}
    };
\end{tikzpicture}\quad
\\
\vspace{0.1in}

\begin{tikzpicture}[thick,scale=0.25]
    \draw (120:3) -- (240:3);
    \draw (120:3) -- (210:3);
    \draw (0:3) -- (30:3);
    \draw (210:3) -- (240:3);
  \draw [black, fill=blue, opacity=0.25] (210:3) -- (240:3) -- (120:3) -- cycle;
    \draw (90:3) -- (120:3);
    \draw (90:3) -- (210:3);
    \draw [black, fill=blue, opacity=0.25] (210:3) -- (90:3) -- (120:3) -- cycle;
    \draw (210:3) -- (330:3);
    \draw (330:3) -- (30:3);
    \draw (270:3) -- (300:3);
    \draw[very thick] (120:3) -- (270:3);
    \draw \foreach \x in {0,30,...,330}
    {
      (\x:3)  node {}
    };
\end{tikzpicture}\quad
\begin{tikzpicture}[thick,scale=0.25]
    \draw (120:3) -- (240:3);
    \draw (120:3) -- (210:3);
    \draw (0:3) -- (30:3);
    \draw (210:3) -- (240:3);
  \draw [black, fill=blue, opacity=0.25] (210:3) -- (240:3) -- (120:3) -- cycle;
    \draw (90:3) -- (120:3);
    \draw (90:3) -- (210:3);
    \draw [black, fill=blue, opacity=0.25] (210:3) -- (90:3) -- (120:3) -- cycle;
    \draw (210:3) -- (330:3);
    \draw (330:3) -- (30:3);
    \draw (270:3) -- (300:3);
    \draw (120:3) -- (270:3);
    \draw[very thick] (30:3) -- (270:3);
    \draw \foreach \x in {0,30,...,330}
    {
      (\x:3)  node {}
    };
\end{tikzpicture}\quad
\begin{tikzpicture}[thick,scale=0.25]
    \draw (120:3) -- (240:3);
    \draw (120:3) -- (210:3);
    \draw (0:3) -- (30:3);
    \draw (210:3) -- (240:3);
  \draw [black, fill=blue, opacity=0.25] (210:3) -- (240:3) -- (120:3) -- cycle;
    \draw (90:3) -- (120:3);
    \draw (90:3) -- (210:3);
    \draw [black, fill=blue, opacity=0.25] (210:3) -- (90:3) -- (120:3) -- cycle;
    \draw (210:3) -- (330:3);
    \draw (330:3) -- (30:3);
    \draw (270:3) -- (300:3);
    \draw (120:3) -- (270:3);
    \draw (30:3) -- (270:3);
    \draw[very thick] (90:3) -- (300:3);
    \draw \foreach \x in {0,30,...,330}
    {
      (\x:3)  node {}
    };
\end{tikzpicture}\quad
\begin{tikzpicture}[thick,scale=0.25]
    \draw (120:3) -- (240:3);
    \draw (120:3) -- (210:3);
    \draw (0:3) -- (30:3);
    \draw (210:3) -- (240:3);
  \draw [black, fill=blue, opacity=0.25] (210:3) -- (240:3) -- (120:3) -- cycle;
    \draw (90:3) -- (120:3);
    \draw (90:3) -- (210:3);
    \draw [black, fill=blue, opacity=0.25] (210:3) -- (90:3) -- (120:3) -- cycle;
    \draw (210:3) -- (330:3);
    \draw (330:3) -- (30:3);
    \draw (270:3) -- (300:3);
    \draw (120:3) -- (270:3);
    \draw (30:3) -- (270:3);
    \draw (90:3) -- (300:3);
    \draw[very thick] (90:3) -- (30:3);
    \draw \foreach \x in {0,30,...,330}
    {
      (\x:3)  node {}
    };
\end{tikzpicture}\quad
\begin{tikzpicture}[thick,scale=0.25]
    \draw (120:3) -- (240:3);
    \draw (120:3) -- (210:3);
    \draw (0:3) -- (30:3);
    \draw (210:3) -- (240:3);
  \draw [black, fill=blue, opacity=0.25] (210:3) -- (240:3) -- (120:3) -- cycle;
    \draw (90:3) -- (120:3);
    \draw (90:3) -- (210:3);
    \draw [black, fill=blue, opacity=0.25] (210:3) -- (90:3) -- (120:3) -- cycle;
    \draw (210:3) -- (330:3);
    \draw (330:3) -- (30:3);
    \draw (270:3) -- (300:3);
    \draw (120:3) -- (270:3);
    \draw (30:3) -- (270:3);
    \draw (90:3) -- (300:3);
    \draw (90:3) -- (30:3);
    \draw[very thick] (30:3) -- (120:3);
    \draw [black, fill=blue, opacity=0.25] (30:3) -- (90:3) -- (120:3) -- cycle;
    \draw [black, fill=blue, opacity=0.25] (30:3) -- (270:3) -- (120:3) -- cycle;
    \draw \foreach \x in {0,30,...,330}
    {
      (\x:3)  node {}
    };
\end{tikzpicture}\quad
\caption{The beginning of a random flag complex process.}
\label{fig:flag}
\end{figure}
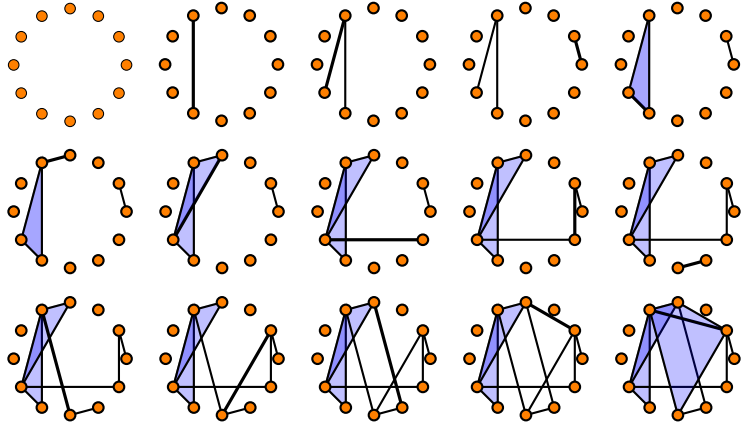

\begin{figure}  
\begin{centering}
\includegraphics[width=3.5in]{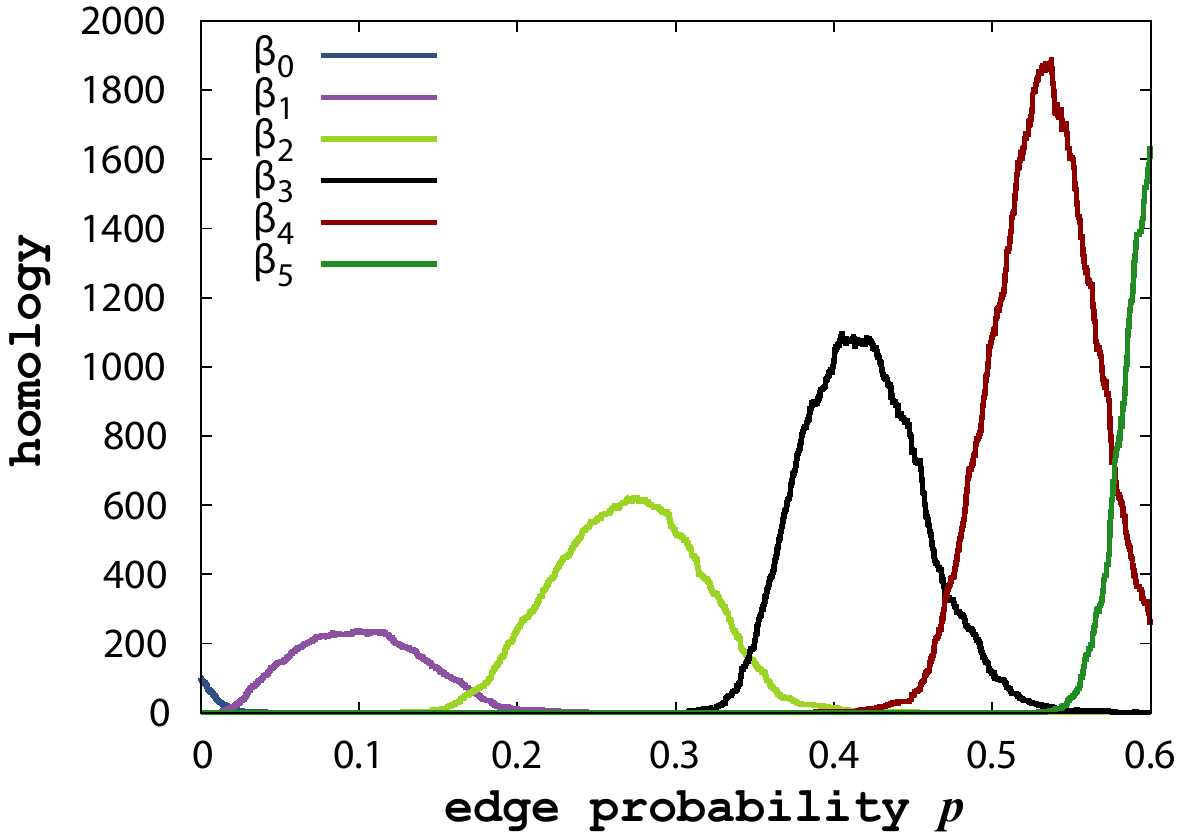}\\
\end{centering}
\caption{The Betti numbers for the beginning of a random flag complex process on $n=100$ vertices. {\it \small (Computation and image courtesy of Afra Zomorodian.)}} 
\label{fig:gnp}
\end{figure}

One can also consider a random flag complex process, which is the same probability space as the random graph process, edges being added one at a time. See Figure \ref{fig:flag}. The Betti numbers of an instance of such a process on $n=100$ vertices and roughly $3000$ steps are illustrated in Figure \ref{fig:gnp}.

We see immediately that homology no longer behaves in a monotone way with respect to the underlying parameter.

\subsection{Vanishing homology}

First of all, we check that, as suggested by Figure \ref{fig:gnp}, there is a range of $p$ outside of which $H_k = 0$ with high probability \cite{clique}.

For the sake of simplicity, we state theorems in the next few sections with assuming $p = n^{-\alpha}$.

\begin{theorem} \label{thm:vanish}
Let $k \ge 1$ and $\alpha > 0$ be fixed, $p= n^{-\alpha}$, and $X \sim X(n,p)$.
\begin{enumerate}
\item If $\alpha > 1/k$ then w.h.p. $H_k(X, \Z)=0$.
\item If $\alpha < 1/(2k + 1)$, then w.h.p. $H_k(X, \Z)=0$.
\end{enumerate}
\end{theorem}

The proof of (1)  is based on showing first that homology is supported on cycles of small support (bounded in size as $n \to \infty$), and then showing that every such cycle is a boundary.

The key observation for (2) is that link of a vertex in a random flag complex is a random flag complex with shifted parameter.  Indeed, even intersecting the links of several vertex links results in another random flag complex.  Then the Nerve Lemma allows one to bootstrap local information about connectivity of a large number of random graphs into global information about cohomology vanishing.  This argument shows something stronger topologically: that if $p=n^{-\alpha}$ where $\alpha < 1/(2k+1)$  then w.h.p.\ $X$ is $k$-connected, i.e.\ $\pi_i(X) = 0$ for $i \le k$.  Recent work of Babson shows that this exponent $1/3$ is tight when $k=1$ \cite{Babson}.

In Section \ref{sec:sharp} we will see that the exponent in (2) can be improved, with a spectral gap argument, if one relaxes to cohomology with rational coefficients.

\subsection{Nontrivial homology and cohomology}

It is also known that for every $k \ge 0$ there is a range of $p=p_k(n)$ for which $H_k(X(n,p)) \neq 0$ with high probability.  Here are three ideas for how one might try to prove this.

\begin{enumerate}
\item {\bf Linear algebra.}  Let $f_i$ denote the number of $i$-dimensional faces of $X$.  Then if $f_k > f_{k-1} + f_{k+1}$, we already have that $H_k \neq 0$ for dimensional reasons, i.e.\ $$\beta_k \ge -f_{k-1} + f_k - f_{k+1}.$$
\item {\bf Homological argument: sphere.}  Try to find a subcomplex $Y$ homeomorphic to a sphere $S^k$, and a simplicial map $f: X \to Y$ such that $f \mid_Y = \mbox{id}$.  Then the homology of $Y$ is naturally a summand of the homology of $X$, and in particular $H_k(X) \neq 0$.
\item {\bf Cohomological argument: isolated face.}  If $\sigma$ is a $k$-dimensional face not contained in any $(k+1)$-dimensional face, then the characteristic function of $\sigma$ represents a cocycle.  If one can somehow show that this function is not also a coboundary, then one has a nontrivial class.
\end{enumerate}

All three of these approaches work, and in roughly the same range of parameter.  They also work equally well for with any choice of coefficients.  Any of these approaches yields the following, for example.

\begin{theorem} \label{thm:nonvanish}
Let $\alpha > 0$ be fixed, $p = n^{-\alpha}$, and $X \sim X(n,p)$.
If $1/ (k+1) < \alpha < 1/ k$, 
then w.h.p.\ $H_k \neq 0$.
\end{theorem}

The first two approaches are discussed in \cite{clique}, and the third approach in \cite{flag}.  In Section \ref{sec:sharp} we will see that the third approach has a slight edge on the other two approaches at the upper end of the window of nontrivial homology.  In this case, a much sharper estimate may be obtained.  So the exponent $1/k$ in Theorems \ref{thm:vanish} and \ref{thm:nonvanish} is sharp.  The exponent $1/(k+1)$ in Theorem \ref{thm:nonvanish} is also sharp, as we will see in Section \ref{sec:sharp}.

\subsection{Limit theorems}  Much more can be shown in the nontrivial regime.  Not only do we know that $H^k \neq 0$ w.h.p., but we can also understand the expectation of $\beta_k$ and its limiting distribution.

The following asymptotic formula for the expectation follows from the linear algebra approach described above.

\begin{theorem}
Let $\alpha > 0$ be fixed, $p = n^{-\alpha}$, and $X \sim X(n,p)$. If $ 1/(k+1) < \alpha < 1/k,$ 
then $$\frac{\mathbb{E}[\beta_k]}{{n \choose k+1} p^{k+1 \choose 2}} \to 1,$$
as $n \to \infty$.
\end{theorem}

Here $\mathbb{E}[\beta_k]$ denotes the expectation of $\beta_k$.  A similar formula can be given for the asymptotic variance $\var[\beta_k]$.

The following central limit theorem characterizes the limiting distribution \cite{Meckes}.

\begin{theorem}[Kahle--Meckes]
Let $\alpha > 0$ be fixed, $p = n^{-\alpha}$, and $X \sim X(n,p)$. If $ 1/(k+1) < \alpha < 1/k$,
then $$\frac{\beta_k - \mathbb{E}[\beta_k]}{\sqrt{ \mbox{Var}[\beta_k]}  } \to \mathcal{N}(0,1)$$
as $n \to \infty$.
\end{theorem} 

Here $ \mathcal{N}(0,1)$ is the standard normal distribution with mean $0$ and variance $1$, and the convergence is in distribution.

\subsection{Sharp thresholds for rational cohomology} \label{sec:sharp}

The following gives a sharp upper threshold for rational cohomology \cite{flag} of random flag complexes.
The Erd\H{o}s-Renyi theorem corresponds to the $k=0$ case.

\begin{theorem} \label{thm:prelim} Let $k \ge 1$ and $\epsilon > 0$ be fixed, and $X = X(n,p)$. 
\begin{enumerate}
\item If $$p \ge  \left( \frac{(k/2 + 1+\epsilon)\log{n}}{n} \right)^{1/ (k+1)},$$
 then $$\prob [ H^k(X, \Q)  =0] \to 1,$$
 \item and if  $$ n^{-1/k + \epsilon} \le p \le  \left( \frac{(k/2 + 1-\epsilon)\log{n}}{n}\right)^{1/(k+1)},$$ then $$\prob [ H^k(X, \Q)  =0] \to 0,$$ 
as $n \to \infty.$
\end{enumerate}
\end{theorem}

The main tools used to prove Theorem \ref{thm:prelim} this are Theorem \ref{thm:ergap1} above which gives a concentration result for the spectral gap, together with the following Theorem \ref{thm:BS}.

\begin{theorem} [Garland, Ballman--\'Swi\k{a}tkowski] \label{thm:BS}  Let $\Delta$ be a pure $D$-dimensional finite simplicial complex such that for every $(D-2)$-dimensional face $\sigma$, the link $\lk_{\Delta}(\sigma)$ is connected and has spectral gap is at least $\lambda_2[ lk_{\Delta}(\sigma)] > 1 - 1 / D$.  Then $H^{D-1}(\Delta,\Q) = 0$.
\end{theorem}

Theorem \ref{thm:BS} is a special case of Theorem 2.5 of Ballman--\'Swi\k{a}tkowski \cite{Ballmann}, which in turn based on earlier work of Garland.  For a deeper discussion of Garland's method, see A.\ Borel's account in S{\'e}minaire Bourbaki \cite{Borel75}.  It is worth noting that Kazhdan had already proved certain cases of Serre's conjecture in 1967 \cite{Kazhdan}, and that this is also the paper in which he introduced Property (T).   


As a corollary, many random flag complexes have nontrivial rational homology only in middle degree.

\begin{corollary} {\label{cor:ddim}} 
Let $d \ge 1$ and $1/(d+1) <  \alpha < 1/d$ be fixed. If $p = n^{-\alpha}$ and $X \sim X(n,p)$, 
then 
$$\homol_i (X, \Q)  = 0 \text{ unless } i =d $$
\end{corollary}

%

\bigskip

It is conceivable that Theorem \ref{thm:prelim} could be sharpened to the following.

\begin{conjecture}
If $$ p = \left( \frac{ (k/2+1) \log{n} + (k/2) \log\log{n} + c}{n} \right)^{1 / (k+1)},$$ where $c \in \R$ is constant, then the dimension of $k$th cohomology $\beta_k$ approaches a Poisson distribution with mean
$$\mu =  \frac{ (k/2+1)^{k/2}}{(k+1)!} e^{-c}.$$  
In particular, $$\prob [ H^k (X, \Q) = 0 ] \to \exp \left[ - \frac{ (k/2+1)^{k/2}}{(k+1)!} e^{-c}\right],$$
as $n \to \infty$.
\end{conjecture}

This conjecture is equivalent to showing that in this regime, cohomology is w.h.p.\ generated by characteristic functions on isolated $k$-faces. This is analogous to the fact that when $c \in \R$ is constant, and $G \sim G(n,p)$ where 
$$p = \frac{ \log n + c}{ n},$$
w.h.p.\ $G$ consists of a giant component and isolated vertices.

\subsection{Torsion}

The question of torsion in random homology is still fairly mysterious, but for certain models of random simplicial complex, we know that torsion will be quite large. It may be surprising to learn, for example, that there exists a $2$-dimensional $\Q$-acyclic simplicial complex $S$ on $31$ vertices with $H_1(S, \Z) $ cyclic of order $$ | H_1(S, \Z) | = 736712186612810774591.$$

The complex is easy to define. The vertices are the elements of the cyclic group $\Z / 31$, the $1$-skeleton is a complete graph, and a set of three vertices $\{ x, y, z \}$ span a $2$-dimensional face if and only if $$x + y + z \equiv 1, 2, \mbox{or } 9 \, (\mbox{mod } 31).$$
This type of ``sum complex'' was introduced and proved to be $\Q$-acyclic using discrete Fourier analysis by Linial, Meshulam, and Rosenthal \cite{LM10}, and I found this particular example choosing $1, 2, 9$ randomly with a calculation in Sage \cite{Sage}.

Work of Kalai \cite{Gil} showed that for $\Q$-acyclic complexes, the expected size of the torsion group is enormous. For example, for a random $2$-dimensional $\Q$-acyclic complex $S_n$ on $n$ vertices, $$\expect \left[| H_1(S, \Z)| \right] \ge e^{cn^2}$$
for some constant $c>0$. This is in some sense the worst case scenario for torsion: all $2$-dimensional $\Q$-acyclic complexes have 
$$| H_1(S, \Z)| \le e^{Cn^2},$$
for some other constant $C>0$; see for example Proposition 3 in Soul\'e \cite{Soule99}.

\bigskip

Kalai's result tells us about the expected size of these random torsion groups, but thirty years later, not much more seems to be known about their expected structure. One natural conjecture might be that the $p$-parts of these random torsion groups obey approach Cohen--Lenstra distributions in the limit.

The Cohen--Lensta distribution over finite abelian $p$-groups proposes that the probability of every group is inversely proportional to the size of its automorphism group. These measures arise naturally in certain number-theoretic settings and have been the subject of a lot of attention, for example in recent work of Ellenberg, Venkatesh, and Westerland \cite{EVW09}. As Lyons points out \cite{Lyons}, the natural measures on $\Q$-acyclic complexes are certain determinantal measures rather than uniform distributions, so this question should be understood with respect to these measures.

In contrast to random $\Q$-acyclic complexes, one might guess for random $2$-complexes or random flag complexes, would be that there may be a window when $H_i ( X, \bf{k}) = 0$ for every fixed field $\bf{k}$ but such that $H_i(X, \Z) \neq 0$, but that this window is very small. The philosophy is that once homology is a finite group, it should only take a relatively small number of random relations to kill it. In particular I would guess the following.

\begin{conjecture} \label{con:torsion}
Let $d \ge 3$ and $$ \frac{1}{d+1} < \alpha < \frac{1}{d}$$ be fixed. If
$p=n^{-\alpha}$ and $X \sim X(n,p)$, then w.h.p.\ $X$ is homotopy equivalent to a bouquet of $d$-dimensional spheres.
\end{conjecture}

It was shown in \cite{clique} that if 
$\alpha < 1/3,$ then w.h.p.\ $ \pi_1(X)=0$. Since simply-connected Moore spaces are unique up to homotopy equivalence, Conjecture \ref{con:torsion} is equivalent to showing that for this range of $p$, $H_{*}( X)$ is w.h.p.\ torsion-free.

\section{Comments}

Now that several different models of random simplicial complex have been studied, we are starting to see a few common themes emerging.

\subsection{There are at least two kinds of topological phase transitions.} The ``upper'' phase transition where homology or cohomology vanishes seems easier to understand cohomologically.  Examples of this kind of phase transition include the Erd\H{o}s--R\'enyi theorem, the Linial--Meshulam theorem, and Theorem \ref{thm:prelim} above.  These thresholds tend to be sharp, happening in a very narrow window.

The ``lower'' phase transition where homology or cohomology passes from vanishing to nonvanishing seems easier to understand homologically.  Theorem \ref{thm:cycle} characterizes the first appearance of cycles in $G(n,p)$.  The higher-dimensional analogue of this phase transition in random complexes is apparently much more subtle, but interesting recent work studying this phase transition appears in \cite{AL} and \cite{ALLM}.  These thresholds tend to be sharp on one side, not on the other.

For random graphs, the basic philosophy local properties such as ``contains a triangle'' have coarse thresholds, whereas global properties, such as ``connected'', have sharp thresholds. For precise statements, see for example the paper by Friedgut with appendix by Bourgain \cite{FB99}. Since isolated faces generate cohomology near the upper phase transition, by the same argument, this is a global property and we should guess the upper threshold is sharp. On the other hand, homology seems to first appear, supported on small spheres, so the lower threshold may tend to be coarse.

\subsection{Homology and cohomology try to be as small as possible.}   With this motto we mean something more geometric, namely that near the lower threshold, homology tends to be supported on small classes.  For example, homology of random geometric complexes in the sparse regime has a basis of vertex-minimal spheres \cite{geometric}.  This happens for both \v{C}ech and Vietoris--Rips complexes, even though the minimal spheres are combinatorially different in the two cases.  This already accounts for the difference in the formulas for expectation of the Betti numbers in the sparse regime.

On the other hand, near the upper threshold cohomology is generated by small classes, namely characteristic functions on isolated $k$-faces.

The existence of these two phase transitions already implies something about what a theory of random persistent homology will have to look like. For example, after a certain point in a random filtration, death times of whatever cycles remain should be well approximated by an appropriate Poisson process.

\subsection {Nature abhors homology.} (With apologies to Aristotle.) Homology is, after all said to measure the number of ``holes'' in a topological space, and holes are made out of vacuum.  More seriously, it is often the case that unless there is a good reason random homology is forced to be there, then it is likely to vanish or be ``small.''

Suppose you have some measure on ``pairs of random linear maps $f: A \to B$ and $g:B \to C$ satisfying $g \circ f = 0$''.  What can we say about the resulting distribution on $H_B$?  Sometimes you can guess the answer from random linear algebra.

For example, one might guess that if $$\dim A \ll \dim B \ll \dim C$$ then there is a good chance that the map $g: B \to C$ is injective, in which case $H_B = 0$.  Or else perhaps $\ker g$ is merely small, but then this still bounds the size of homology.  Similarly, if $$\dim A \gg \dim B \gg \dim C,$$ then one might expect $f: A \to B$ is probably surjective (or nearly so) and so one expects that $H_B$ is small.  In fact the only place the one place where we actually expect to see large homology, if dimension were the only consideration, would be if $$\dim A \ll \dim B \gg \dim C.$$
So you might guess that $\dim H_B \approx \max \{ 0,  - \dim A + \dim B - \dim C \}$.

The random flag complex $X(n,p)$ is an example of where this kind of argument works surprisingly well. In a random flag complex, the theorems so far tell us that we should only expect to see one or two nontrivial Betti numbers at any given part in the process, and that the overlap between two nonzero Betti numbers should not be too large. If we make the simplifying assumption that all of the reduced Betti numbers are zero except for one, then the  nonzero Betti number equals the absolute value of the reduced Euler characteristic. The upshot is that the expected Euler characteristic $\expect[ \tilde{\chi} ]$ is easy to compute, as linearity of expectation gives that 
$$ \expect[ \tilde{\chi} ] = - 1 + n - {n \choose 2} p + { n \choose 3} p^3 - {n \choose 4} p^6 + {n \choose 5} p^{10} - \dots.$$

\begin{figure} 
\begin{centering}
\includegraphics[width=3.5in]{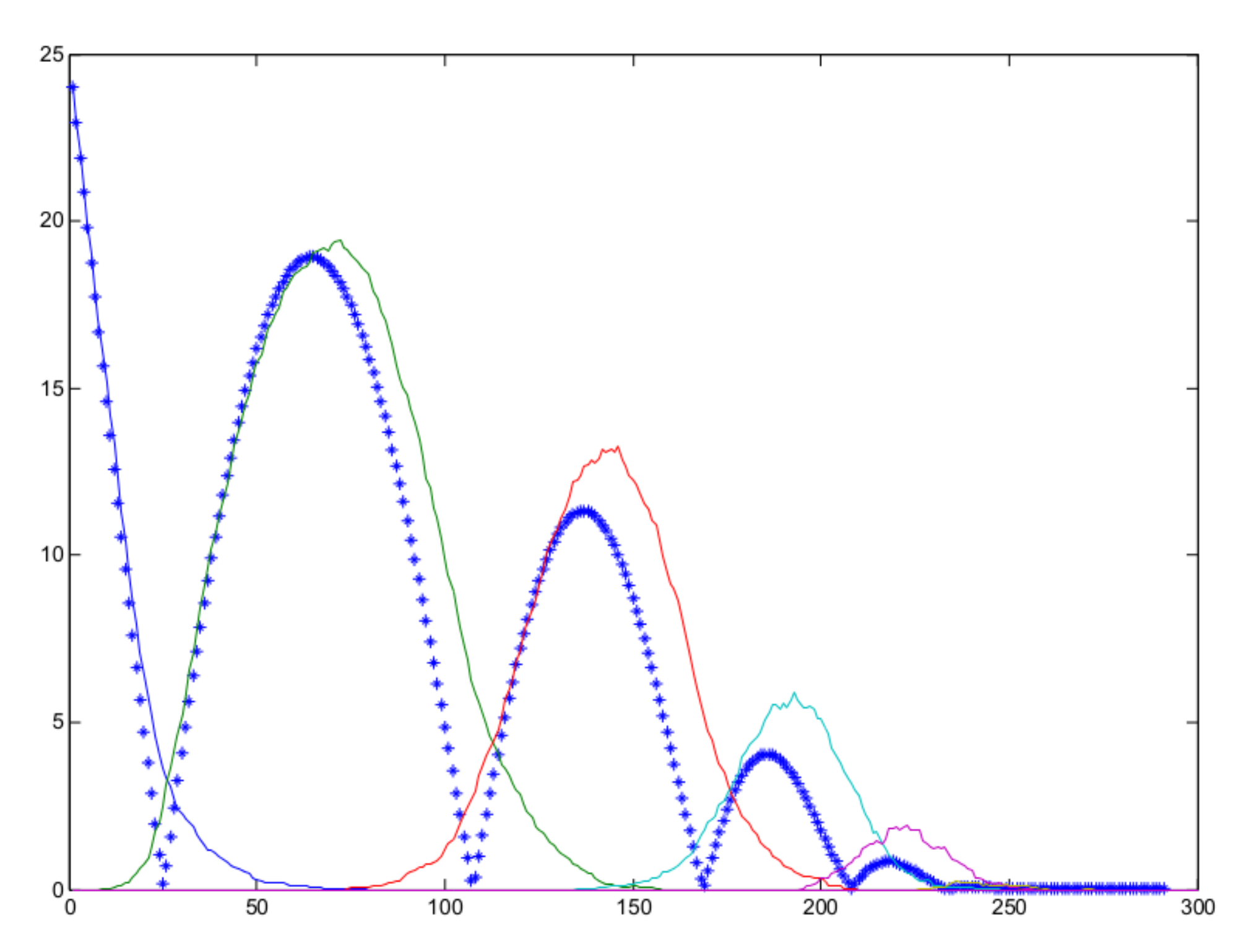}\\
\end{centering}
\caption{The average of $50$ complete random flag complex processes on $n=25$ vertices (continuous curves), plotted against $|\expect[ \tilde{\chi} ]|$ (blue stars). As in Figure \ref{fig:gnp}, the horizontal axis is the number of edges and the vertical axis the Betti numbers. {\it \small (Computation and image courtesy of Vidit Nanda.)}} 
 \label{fig:n25}
\end{figure}

See, for example, Figure \ref{fig:n25}. Using his software Perseus \cite{perseus}, Vidit Nanda computed fifty complete flag complex processes on $n=25$ vertices, and the average is plotted against the expected Euler characteristic. It is quite striking that the prediction is so good, especially considering that all the theory is as $n \to \infty$ and $p \to 0$, but here the prediction seems to work reasonably well, even for relatively small $n$ and large $p$.

\subsection{Random simplicial complexes are expanders}  Random simplicial complexes have many expander-like properties.
Higher-dimensional analogues of expanders have attracted a lot of attention recently--- see for example the discussion in \cite{DK11}, and also Gromov's recent work on ``geometric overlap'' properties of expanders  \cite{Gromov1, Gromov2}.

For a geometric application, expander graphs are well understood to be impossible to embed in Euclidean space with small metric distortion, by work of Bourgain \cite{Bourgain85}. Volume distortion analogues of this were recently established for random simplicial complexes by Newman--Rabinovich \cite{volume} and by Dotterrer \cite{Dominic}.
 
We would like to better understand the higher-dimensional analogues of the Cheeger--Buser inequalities relating spectral gap and the ``bottleneck'' expansion constant.  For recent work on higher-dimensional analogues of Cheeger--Buser, see Gundert--Wagner \cite{GW12}, and Steenbergen et al.\ \cite{SKM}.

\subsection{A multi-parameter model.}  A model that deserves more attention is the multi-parameter random complex $\Delta(n; p_1, p_2, \dots)$.  Here there are $n$ vertices, the probability of an edge is $p_1 = p_1(n)$, and the complex is built inductively by dimension in so that the probability of every $k$-dimensional simplex, conditioned that its entire $(k-1)$-dimensional boundary is already in place, is $p_k= p_k(n)$, independently.

Several of the random simplicial complexes discussed here are special cases of this model; in particular, the random graph 
$$G(n,p)=\Delta(n; p_1, 0, 0, \dots),$$ 
the random $2$-complex
$$Y(n,p) = \Delta(n; 1, p_2, 0, 0, \dots),$$
and the random flag complex
$$X(n,p)=\Delta(n; p_1, 1, 1,  \dots).$$ 

\section*{Acknowledgements} I am grateful to IAS for hosting me for two weeks in January 2013, and for having a chance to present this material in two talks at the MacPherson seminar while the article was in preparation, and to receive feedback from a friendly audience. I also thank Peter Landweber for a careful read of an earlier draft, and for several helpful suggestions. 

\bigskip

I dedicate this article to Gunnar Carlsson, and thank him for his encouragement and support.

\bigskip

\bibliographystyle{amsplain}
\bibliography{specrefs}


\end{document}